\documentclass{amsart}

\usepackage{a4,xy,amssymb,diagrams,times,array,color}
\xyoption{poly}
\xyoption{arc}
\xyoption{knot}
\xyoption{curve}
\xyoption{arrow}
\xyoption{all}

\makeindex

\newcommand{\C}{\mathbb{\Bbbk}}

\newcommand{\wis}[1]{{\text{\em \usefont{OT1}{cmtt}{m}{n} #1}}}
\newcommand{\Oscr}{\mathcal{O}}

\newcommand{\vtx}[1]{*+[o][F-]{\scriptscriptstyle #1}}

\newcolumntype{C}{>{$}c<{$}}

\newcommand{\aar}[2]{ \ar@{-}|(.25){{\bf #1}}|{\bullet}|(.75){{\bf  #2}}}
\newtheorem{definition}{Definition}

\newtheorem{theorem}{Theorem}

\title{Noncommutative Geometry and Dual Coalgebras}
\author{Lieven Le Bruyn} 
\address{Department of Mathematics, University of Antwerp \\ 
 Middelheimlaan 1, B-2020 Antwerp (Belgium) \\ {\tt lieven.lebruyn@ua.ac.be}}

\begin{document}
\sloppy
 
\def\ldb{\mathopen{\{\!\!\{}} \def\rdb{\mathclose{\}\!\!\}}}

 %\mathversion{bold}

 \begin{abstract}
In arXiv:math/0606241v2 M. Kontsevich and Y. Soibelman argue that the category of noncommutative (thin) schemes is equivalent to the category of coalgebras. We propose that under this correspondence the affine scheme $\wis{rep}(A)$ of a $k$-algebra $A$ is the dual coalgebra $A^o$ and draw some consequences. In particular, we describe the dual coalgebra $A^o$ of $A$ in terms of the $A_{\infty}$-structure on the Yoneda-space of all the simple finite dimensional $A$-representations.
\end{abstract}

\maketitle

\tableofcontents

\section{$\wis{rep}(A)=A^o$}

Throughout, $k$ will be a (commutative) field with separable closure $\C$. In \cite[\S I.2]{KontSoib} Maxim Kontsevich and Yan Soibelman define a {\em noncommutative thin scheme} to be a covariant functor commuting with finite projective limits
\[
\wis{X}~:~\wis{alg}_k^{fd} \rTo \wis{sets} \]
from the category $\wis{alg}_k^{fd}$ of all {\em finite dimensional} $k$-algebras (associative with unit) to the category $\wis{sets}$ of all sets. They prove \cite[Thm. 2.1.1]{KontSoib} that every noncommutative thin scheme is represented by a $k$-coalgebra. 

Recall that a $k$-coalgebra is a $k$-vectorspace $C$ equipped with linear structural morphisms :  a comultiplication $\Delta~:~C \rTo C \otimes C$ and a counit $\epsilon~:~C \rTo k$ satisfying the coassociativity $(id \otimes \Delta) \Delta = (\Delta \otimes id) \Delta$ and counitary property $(id \otimes \epsilon) \Delta = (\epsilon \otimes id) \Delta = id$. 

By being representable they mean that every noncommutative thin scheme $\wis{X}$ has associated to it a $k$-coalgebra $C_{\wis{X}}$ with the property that for any finite dimensional $k$-algebra $B$ there is a natural one-to-one correspondence
\[
\wis{X}(B) = \wis{alg}_k(B,C_{\wis{X}}^*) \]
Here, for a $k$-coalgebra $C$ we denote by $C^*$ the space of linear functionals $Hom_k(C,k)$ which acquires a $k$-algebra structure by dualizing the structural coalgebra morphisms. 

They call $C_{\wis{X}}$ the {\em coalgebra of distributions} on $\wis{X}$ and define the {\em noncommutative algebra of functions} on $\wis{X}$ to be the dual $k$-algebra $k[\wis{X}] = C_{\wis{X}}^*$.

Whereas the dual $C^*$ of a $k$-coalgebra is always a $k$-algebra, for a $k$-algebra $A$ it is not true in general that the dual vectorspace $A^*$ is a coalgebra, because $(A \otimes A)^* \not\simeq A^* \otimes A^*$. Still, one can define the subspace
\[
A^o = \{ f \in A^* = Hom_k(A,k)~|~ker(f)~\text{contains a twosided ideal of finite codimension}~\} \]
and show that the duals of the structural morphisms on $A$ determine a $k$-coalgebra structure on this dual coalgebra $A^o$, see for example \cite[Prop. 6.0.2]{Sweedler}. 

With these definitions, {\em Kostant duality} asserts that the functors
\[
\xymatrix{\wis{alg}_k \ar@/^/[rr]^{o} & & \wis{coalg}_k \ar@/^/[ll]^{\ast}}
\]
are adjoint, \cite[Thm. 6.0.5]{Sweedler}. That is, for any $k$-algebra $A$ and any $k$-coalgebra $C$, there is a natural one-to-one correspondence between the homomorphisms
\[
\wis{alg}_k(A,C^*) = \wis{coalg}_k(C,A^o) \]
Moreover, we have \cite[Lemma 6.0.1]{Sweedler} that for $f \in \wis{alg}_k(A,B)$, the dual map $f^*$ determines a $k$-coalgebra morphism $f^* \in \wis{coalg}_k(B^o,A^o)$.

For a $k$-algebra $A$ one can define the contravariant functor $\wis{rep}(A)$ describing its finite dimensional representations \cite[Example 2.1.9]{KontSoib}
\[
\wis{rep}(A)~:~\wis{coalg}_k^{fd} \rTo \wis{sets} \qquad C \mapsto \wis{alg}_k(A,C^*) \]
from {\em finite dimensional} $k$-coalgebras $\wis{coalg}_k^{fd}$ to $\wis{sets}$, which commutes with finite direct limits. As on finite dimensional $k$-(co)algebras  Kostant duality is an anti-equivalence of categories
\[
\xymatrix{\wis{alg}_k^{fd} \ar@/^/[rr]^{\ast} & & \wis{coalg}_k^{fd} \ar@/^/[ll]^{\ast}}
\]
we might as well describe $\wis{rep}(A)$ as the noncommutative thin scheme represented by $A^o$
\[
\wis{rep}(A)~:~\wis{alg}_k^{fd} \rTo \wis{sets} \qquad B=C^* \mapsto \wis{alg}_k(A,B=C^*) = \wis{coalg}_k(C=B^*,A^{o}) \]
the latter equality follows again from Kostant duality. Therefore, we propose

\begin{definition} The {\em noncommutative affine scheme} $\wis{rep}(A)$ is the noncommutative (thin) scheme corresponding to the dual $k$-coalgebra $A^o$ of $A$.
\end{definition}

\section{$\wis{simp}(A) = corad(A^o)$}

The dual $k$-coalgebra $A^o$ is usually a {\em huge} object and hence contains a lot of information about the $k$-algebra $A$. Let us begin by recalling how the geometry of a commutative affine $k$-scheme $X$ is contained in the dual coalgebra $A^o$ of its coordinate ring $A = \C[X]$. 

Recall that a coalgebra $D$ is said to be {\em simple} if it has no proper nontrivial subcoalgebras. In particular, a simple coalgebra $D$ is finite dimensional over $k$ and by duality is such that $D^*$ is a simple $k$-algebra, that is, $D^*$ is a central simple $L$-algebra where $L$ is a finite separable extension of $k$.

Hence, in case $A=\C[X]$ (and $\C$ is separably closed) we have that all simple subcoalgebras of $A^o$ are one-dimensional (and hence are spanned by a group-like element), because they correspond to simple representations of $A$. 

That is, $A^o$ is {\em pointed} and by \cite[Prop. 8.0.7]{Sweedler} we know that any cocommutative pointed coalgebra is the direct sum of its {\em pointed irreducible components} (at the algebra level, this says that a semi-local commutative algebra is the direct sum of locals). Therefore,
\[
A^o = \oplus_{x \in X} C_x \]
where each $C_x$ is pointed irreducible and cocommutative. As such, each $C_x$ is a subcoalgebra of the enveloping coalgebra of the abelian Lie algebra on the tangent space $T_x(X)$. That is, we recover the points of $X$ as well as tangent information from the dual coalgebra $A^o$.

But then, the dual algebra of $A^o$, that is the 'noncommutative' algebra of functions $A^{o\ast}$  decomposes as
\[
A^{o\ast} = \oplus_{x \in X} \hat{\Oscr}_{x,X} \]
the direct sum of the completions of the local algebras at points. The diagonal embedding $A= \C[X] \rInto A^{o \ast}$ inevitably leads to the structure scheaf $\Oscr_X$.

We will now associate a topological space associated to any $k$-algebra $A$, generalizing the space of points equipped with the Zariski topology when $A$ is a commutative affine $k$-algebra. In the next section we will describe the dual coalgebra $A^o$ when $A$ is a noncommutative affine $\C$-algebra.

The {\em coradical} $corad(C)$ of a $k$-coalgebra $C$ is the (direct) sum of all simple subcoalgebras of $C$. It is also the direct sum of all simple subcomodules of $C$, when $C$ is viewed as a left (or right) $C$-comodule.

In the example above, when $A= \C[X]$, we have that $corad(A^o) = \oplus_{x \in X} \C~ev_x$ where the group-like element $ev_x$ is evaluation in the point $x$. This motivates :

\begin{definition} For a $k$-algebra $A$ we define the {\em space of points} $\wis{simp}(A)$ to be the set of direct summands of $corad(\wis{rep}(A)) = corad(A^o)$. That is, $\wis{simp}(A)$ is the set of simple subcoalgebras of $\wis{rep}(A)$.
\end{definition}

By Kostant duality it follows that $\wis{simp}(A)$ is the set of all finite dimensional simple algebra quotients of the $k$-algebra $A$, or equivalently, the set of all isomorphism classes of finite dimensional simple $A$-representations, explaining the notation.

We can equip this set with a {\em Zariski topology} in the usual way, using the {\em evaluation map}
\[
A^o \times A \rTo^{ev} k \qquad (f,a) \mapsto f(a) \]
when restricted to the subcoalgebra $corad(A^o)$. Note that the evaluation map actually defines a {\em measuring} of $A$ to $k$ \cite[Prop. 7.0.3]{Sweedler}, that is,
$A^o \otimes A \rTo^{ev} k$ 
satisfies
\[
ev(f \otimes aa') = \sum_{(f)} f_{(1)}(a) f_{(2)}(a') \qquad \text{and} \qquad ev(f \otimes 1) = \epsilon(f) 1_k\]

\begin{definition} The {\em Zariski topology} of a $k$-algebra $A$ is the set $\wis{simp}(A)$ equipped with the topology generated by the basic closed sets
\[
\mathbb{V}(a) = \{ S \in \wis{simp}(A)~|~ev(S \otimes a) = 0,~\text{that is}~f(a)=0,~\forall f \in S~\} \]
\end{definition}

Having associated a topological space to a $k$-algebra, one might ask when this is a functor.
Functoriality has always been a problem in noncommutative geometry. Indeed, a simple $B$-representation does not have to remain a simple $A$-representation under restriction of scalars via $\phi~:~A \rTo B$. 

Still, if we define $\wis{rep}(A)=A^o$, we get functionality for free. If $A \rTo^{\phi} B$ is an algebra morphism, we have seen that the dual map maps $B^o$ to $A^o$, so we have a morphism
\[
B^o =  \wis{rep}(B) \rTo^{\phi^*} \wis{rep}(A) = A^o \]
A coalgebra is the direct limit of its finite dimensional coalgebras, and they correspond under duality to finite dimensional algebras. Hence, $\phi^*$ is  the natural map on finite dimensional representations  by restriction of scalars.

The observed failure of functoriality on the level of points translates
on the coalgebra-level to the fact that for a coalgebra map $B^o \rTo A^o$ the coradical $corad(B^o)$ does not have to be mapped to $corad(A^o)$, in general. 

However, when $corad(B^o)$ is cocommutative, we {\em do have} that $\phi^*(corad(B^o)) \subset corad(A^o)$ by \cite[Thm. 9.1.4]{Sweedler}. In particular, we recover the functor of points in {\em commutative} algebraic geometry.

Clearly, we still have  $corad(B^o) \rTo A^o$ in general. This corresponds to the fact that there is always a map $\wis{simp}(B) \rTo \wis{rep}(A)$.

Next, let us turn to the algebra of functions on $\wis{rep}(A)$. By definition we have
\[
k[\wis{rep}(A)] = A^{o \ast} \]
and we can ask how this algebra relates to the algebra $A$.

In general, it is {\em not} true that $A \rInto A^{o\ast}$. This only holds when $A^o$ is dense in $A^*$ in which case the $k$-algebra is said to be {\em proper}, see \cite[\S 6.1]{Sweedler}.

In the commutative case, when $A$ is a finitely generated $k$-algebra, then $A$ is indeed proper and this is a consequence of the Hilbert Nullstellensatz and the Krull intersection theorem. 

When $A$ is noncommutative, this is no longer the case. For example, if $A = A_n(k)$ the {\em Weyl algebra} over a field of characteristic zero $k$, then $A$ is simple whence has no twosided ideals of finite codimension. As a result $A^o = 0$! As our proposal for the noncommutative affine scheme $\wis{rep}(A)$ is based on finite dimensional representations of $A$, it will not be suitable for $k$-algebras having few such representations.

\section{the dual coalgebra $A^o$}

In general though, $A^o$ is a huge object, so it is very difficult to describe explicitly. In this section, we will begin to tame $A^o$ even when $A$ is noncommutative. 

In order not to add extra problems, we will assume that $\C$ is separably closed in this section. The general case can be recovered by taking $Gal(\C/k)$-invariants (replacing quivers by {\em species} in the sequel).

Over a separably closed field $\C$ all simple subcoalgebras are full matrix coalgebras $M_n(\C)^*$, that is, $M_n(\C)^* = \oplus_{i,j} \C e_{ij}$ with $\Delta(e_{ij}) = \sum_{k=1}^n e_{ik} \otimes e_{kj}$ and $\epsilon(e_{ij}) = \delta_{ij}$. 

Hence, $corad(A^o) = \oplus_S M_{n_S}(\C)^*$ where the sum is taken over all finite dimensional simple $A$-representations $S$, each having dimension $n_S$. 

In algebra, one can resize idempotents by Morita-theory and hence replace full matrices by the basefield.  In coalgebra-theory there is an analogous duality known as {\em Takeuchi equivalence}, see \cite{Takeuchi}.

The isotypical decomposition of $corad(A^o)$ as an $A^o$-comodule is of the form $\oplus_S C_S^{\oplus n_S}$, the sum again taken over all simple $A$-representations. Take the $A^o$-comodule
$E= \oplus_S C_S$ and its {\em coendomorphism coalgebra} 
\[
A^{\dagger} = coend^{A^o}(E) \]
then Takeuchi-equivalence (see for example \cite[\S 4, \S 5]{Chin} and the references contained in this paper for more details) asserts that $A^o$ is Takeuchi-equivalent to the coalgebra $A^{\dagger}$ which is {\em pointed}, that is, $corad(A^{\dagger}) = \C~\wis{simp}(A) = \oplus_S \C g_S$ with one {\em group-like} element $g_S$ for every simple $A$-representation. Remains to describe the structure of the full  basic coalgebra $A^{\dagger}$.

For a (possibly infinite) quiver $\vec{Q}$ we define the {\em path coalgebra} $\C \vec{Q}$ to be the vectorspace $\oplus_p \C p$ with basis all oriented paths $p$ in the quiver $\vec{Q}$ (including those of length zero, corresponding to the vertices) and with structural maps induced by
\[
\Delta(p) = \sum_{p=p'p"} p' \otimes p" \qquad \text{and} \qquad \epsilon(p) = \delta_{0,l(p)} \]
where $p'p"$ denotes the concatenation of the oriented paths $p'$ and $p"$ and where $l(p)$ denotes the length of the path $p$. Hence, every vertex $v$ is a group-like element and for an arrow $\xymatrix{\vtx{v} \ar[r]^a & \vtx{w}}$ we have $\Delta(a) = v \otimes a + a \otimes w$ and $\epsilon(a)=0$, that is, arrows are skew-primitive elements.

For every natural number $i$, we define the {\em $ext^i$-quiver} $\overrightarrow{\wis{ext}}^i_A$ to have one vertex $v_S$ for every $S \in \wis{simp}(A)$ and such that the number of arrows from $v_S$ to $v_T$ is equal to the dimension of the space $Ext^i_A(S,T)$. With $\wis{ext}^i_A$ we denote the $\C$-vectorspace on the arrows of $\overrightarrow{\wis{ext}}^i_A$.

The {\em Yoneda-space} $\wis{ext}^{\bullet}_A = \oplus \wis{ext}^i_A$ is endowed with a natural $A_{\infty}$-structure \cite{Keller}, defining a linear map (the {\em homotopy Maurer-Cartan map}, \cite{Segal})
\[
\mu = \oplus_i m_i~:~\C \overrightarrow{\wis{ext}}^1_A \rTo \wis{ext}^2_A \]
from the path coalgebra $\C \overrightarrow{\wis{ext}}^1_A$ of the $ext^1$-quiver to the vectorspace $\wis{ext}^2_A$, see \cite[\S 2.2]{Keller} and \cite{Segal}.

\begin{theorem} The dual coalgebra $A^o$ is Takeuchi-equivalent to the pointed coalgebra $A^{\dagger}$ which is  the sum of all subcoalgebras contained in the kernel of the linear map 
\[
\mu = \oplus_i m_i~:~\C \overrightarrow{\wis{ext}}^1_A \rTo \wis{ext}^2_A \]
 determined by the $A_{\infty}$-structure on the Yoneda-space $\wis{ext}^{\bullet}_A$.
\end{theorem}

We can reduce to finite subquivers as any subcoalgebra is the limit of finite dimensional subcoalgebras and because any finite dimensional $A$-representation involves only finitely many simples. Hence, the statement is a global version of the result on finite dimensional algebras due to B. Keller \cite[\S 2.2]{Keller}. 

Alternatively, we can use the results of E. Segal \cite{Segal}. Let $S_1,\hdots,S_r$ be distinct simple finite dimensional $A$-representations and consider the semi-simple module $M=S_1 \oplus \hdots \oplus S_r$ which determines an algebra epimorphism
\[
\pi_M~:~A \rOnto M_{n_1}(\C) \oplus \hdots \oplus M_{n_r}(\C) = B \]
If $\mathfrak{m} = Ker(\pi_M)$, then the $\mathfrak{m}$-adic completion $\hat{A}_{\mathfrak{m}} = \underset{\leftarrow}{lim}~A/\mathfrak{m}^n$ is an augmented $B$-algebra and we are done if we can describe its finite dimensional (nilpotent) representations. Again, consider the $A_{\infty}$-structure on the Yoneda-algebra $Ext^{\bullet}_A(M,M)$ and the quiver on $r$-vertices $\overrightarrow{\wis{ext}}^1_A(M,M)$ and the homotopy Mauer-Cartan map
\[
\mu_M = \oplus_i m_i~:~\C \overrightarrow{\wis{ext}}^1_A(M,M) \rTo Ext^2_A(M,M) \]
Dualizing we get a subspace $Im(\mu_M^*)$ in the path-{\em algebra} $\C \overrightarrow{\wis{ext}}^1_A(M,M)^*$ of the dual quiver. Ed Segal's main result \cite[Thm 2.12]{Segal} now asserts that $\hat{A}_{\mathfrak{m}}$ is Morita-equivalent to
\[
\hat{A}_{\mathfrak{m}} \underset{M}{\sim} \frac{(\C \overrightarrow{\wis{ext}}^1_A(M,M)^*)^{\hat{}}}{(Im(\mu_M^*))} \]
where $(\C \overrightarrow{\wis{ext}}^1_A(M,M)^*)^{\hat{}}$ is the completion of the path-algebra at the ideals generated by the paths of positive length. The statement above is the dual coalgebra version of this.

\end{document}